# Adaptive Certainty-Equivalence Control With Regulation-Triggered Finite-Time Least-Squares Identification, Part I: Design

**Iasson Karafyllis* and Miroslav Krstic***

*Dept. of Mathematics, National Technical University of Athens,
Zografou Campus, 15780, Athens, Greece,
email: iasonkar@central.ntua.gr

**Dept. of Mechanical and Aerospace Eng., University of California, San Diego, La Jolla, CA 92093-0411, U.S.A., email: krstic@ucsd.edu

**Abstract**

For general nonlinear control systems we present a novel approach to adaptive control, which employs a certainty equivalence (indirect) control law and an identifier with event-triggered updates of the plant parameter estimates, where the triggers are based on the size of the plant's state and the updates are conducted using a non-recursive least-squares estimation over certain finite time intervals, with updates employing delayed measurements of the state. With a suitable non-restrictive parameter-observability assumption, our adaptive controller guarantees global stability, regulation of the plant state, and our identifier achieves parameter convergence, in finite time, even in the absence of persistent excitation, for all initial conditions other than those where the initial plant state is zero. The robustness of our event-triggered adaptive control scheme to vanishing and non-vanishing disturbances is verified in simulations with the assistance of a dead zone-like modification of the update law. The major distinctions of our approach from supervisory adaptive schemes is that our approach is indirect and our triggering is related to the control objective (the regulation error). The major distinction from the classical indirect Lyapunov adaptive schemes based on tuning related to the regulation error is that our approach does not involve a complex redesign of the controller to compensate for the detrimental effects of rapid tuning on the transients by incorporating the update law into the control law. Instead, our approach allows for the first time to use a simple certainty equivalence adaptive controller for general nonlinear systems. All proofs are given in a companion paper.

## 1. Introduction

*New Design Approach.* Classical tuning-based adaptive control went through several spurts of methodological advances from the late 1950s onward, remaining focused on linear systems until about 1990 [10,41,45,20], and culminating in the designs for nonlinear systems in the mid-1990s [27]. Refinements and extensions of classes of systems have followed but paradigm-changing *methodological* advances—namely, fundamentally new *approaches* to adaptive control—have been nearly non-existent in the last two decades.

With this paper we introduce a design approach that is new at a methodological level. In the previous Lyapunov-based and modular (estimation-based) approaches to adaptive control, either the certainty equivalence approach is abandoned by designing complex and interdependent controllers and parameter estimators [27,28,30] or the disturbing effect of parameter estimation is counteracted by slowing down the estimation [10,41,45,20] or through controller strengthening [27,29,31]. In our new design paradigm we employ unconventional *regulation-triggered* identifiers that allow us to use simple certainty equivalence controllers without imposing growth conditions on nonlinearities and without normalizing/slowing adaptation. Crucial to our design are three ideas: (2.1) employing piecewise-constant parameter estimates in the controller between the event-based triggers for the identifier, which eliminates the worry about the disturbing effect of rapidly changing estimates between the estimate switches, (2.2) employing finite-time (dead-beat) parameter identifiers, which eliminate the problems caused by long-lasting large parameter estimation errors, and (3.1) using regulation error to trigger the parameter estimate switches, rather than estimation error-based triggers, which makes the parameter updating rate as rapid as necessary to prevent instability regardless of the growth rates of the plant's nonlinearities.

A fourth feature of our approach, which is not essential but is valuable, is that our finite-time identifier design is based on a least-squares (LS) approach, enabling a balanced convergence rate across the entries of the parameter estimation vector, which is the principal advantage of LS over all other estimator approaches. In other approaches (gradient, Lyapunov, etc) it is impossible to guarantee a priori a balanced convergence because the excitation levels for individual columns of the regressor matrix are not known a priori.

*Existing Adaptive Control Approaches.* To understand the merits (or shortfalls) of any adaptive control design claimed to be methodologically novel, it is crucial to understand the categorizations of adaptive control approaches. Classical and comprehensive references such as [43,28,20] are helpful for this purpose.

First, adaptive controllers can be categorized into indirect (updating plant parameters) and direct (updating controller parameters). The indirect approaches have the advantages of not requiring overparameterization, being more readily applicable to nonlinear systems, and yielding estimates of physical parameters, whereas the direct approaches avoid the on-line solution of "design equations" (such as Bezout, Riccati, etc.).

Second, adaptive controllers are categorized into certainty equivalence (CE) controllers, namely those designed for the "known parameter'' situation, and those that are not of the CE type.

Third, non-CE adaptive controllers generally come in two classes. One class explicitly incorporates the functional form of the parameter update law, in order to compensate the rapidly time varying character of the parameter estimates (such adaptive controllers are usually designed using complete Lyapunov functions and are referred to as "Lyapunov-based" [27,30]). The other class allows off-the-shelf gradient and LS parameter estimators and does not incorporate the update law into the control law but, in order to compensate for the disturbing effects of the parameter estimation error and its rate of variation, which affects the validity of the control law, employs strengthening of the control law against such disturbances in some form of nonlinear damping (such adaptive controllers are referred to as "modular" [29,31] because the controller and the parameter estimator are designed and analyzed separately, or "estimation-based" [43] because the estimator is not based on the regulation or stabilization objective but merely on the control-unrelated objective of estimating the parameters).

All the above approaches to adaptive control are developed for general parabolic Partial Differential Equations in one spatial dimension in the book [46], following the systematization laid out in Chapter 7 of that book.

*Event-Triggered Control.* In our approach we employ event-triggered identifiers. Event-triggered control has attracted considerable attention within the control systems community. Indeed, event-triggered control has provided solutions for difficult control problems that involve sampling, quantized measurements, output-feedback control, distributed networked control and decentralized

control; see [3,5,6,8,9,11,13,14,15,32,33,34,48,49,51,52,53,56]. In all the cases, the system under event-triggered control becomes a hybrid dynamical system.

*Supervisory Approach to Adaptive Control.* Event-triggered forms of adaptive control have existed for over two decades now. There exist direct *adaptive* control schemes in the literature with guaranteed convergence properties for the closed-loop system. Direct adaptive control approaches for linear systems have been proposed in [35,36,37,38]: the proposed schemes either employ event-triggering or sampled-data techniques. Direct adaptive control design methodologies with logic-based switching for linear and nonlinear control systems have been developed in [16,17,18,19,39,40] (see also the references therein): the proposed supervisory adaptive control schemes employ multi-model based estimators of the performance of the "current" controller in conjunction with hierarchical hysteresis switching logic (which is the event-triggered element in the design). This approach is akin to *estimation error-triggered* controller scheduling (in a *direct* adaptive setting), as opposed to our *regulation-triggered* identification (in an *indirect* setting). We estimate the plant's unknown parameters with a dead-beat, least-squares identifier with delays, which allows us to derive (constructively) appropriate KL estimates and employ a simple, single CE control law. Moreover, we are able to guarantee convergence of the parameters estimates in a *pre-specified* time horizon. The advantage of logic-based switching [16,17,18,19] is that it can also deal with systems with disturbances while the results of the present paper are for undisturbed systems.

The robustness of our event-triggered adaptive control approach with respect to vanishing and non-vanishing perturbation is tested by means of a simulation study (see Section 7) and its performance is compared with the achieved performance of conventional adaptive controllers (designed by means of methodologies provided in [20,28]). The results show that the robustness properties of our regulation-triggered adaptive controller are comparable to the robustness properties of the closed-loop system with the nominal controller and known parameter values.

*Least-Squares Identifiers in Adaptive Control.* Least-squares (LS) identification is attractive because of the ability of the Riccati equation to adjust the adaptation gain in real time to the actual signal content of the regressor matrix so that all the channels of the parameter estimation vector evolve at comparable speed, rather than at vastly different time scale (with the slowest one being dominant). Recent advances and applications in LS-based adaptive control are [4,26,42,54,55]. However, in most LS adaptive control methods the parameter convergence is not guaranteed without persistency of excitation (PE), which is seldom verifiable a priori. A different LS estimator, which uses a hybrid dead-beat observer with delays was proposed in [21,23]. The estimator in [21] does not require PE but only the weaker assumption of strong observability and was shown to be robust with respect to measurement noise.

*Our Approach and Contribution of the Paper.* Our adaptive control approach uses the idea of estimating the unknown parameters by means of a *dead-beat* identifier in conjunction with a certainty equivalence controller. As a result, the closed-loop convergence properties ultimately become (after finite time) those of the nominal (known-parameter) controller. To this end, we use the hybrid dead-beat observer proposed in [21] (but slightly modified). The rate of adaptation is determined by a scheme triggered by the regulation error: when things "do not go well" (with the plant's state), the adaptation is accelerated. It is in this manner that the finite escape phenomenon is avoided. Furthermore, all excitation assumptions used in adaptive control literature are replaced by an appropriate observability assumption (see Assumption (H3) below). Our parameter-observability assumption can be verified a priori (see Section 5 as well as the companion paper [24]).

Dead-beat identifiers were first proposed in [1,2,12]. They require a PE assumption and are not of an LS type. The identifiers in in [1,2,12] could conceivably used with our regulation-triggered approach and CE control, but we do not prove such a result here.

*Organization of the Paper.* Section 2 is devoted to the formulation of the problem and the presentation of the assumptions under which the adaptive regulator is constructed. Section 3 provides the detailed description of the event-triggered identifier and the adaptive controller. The main results of the present work are given in Section 4 (Theorem 4.1, Theorem 4.2 and Theorem 4.3 and Corollary 4.4). Section 5 provides a simple nonlinear example that illustrates the design of the adaptive controller as well as the verification of all assumptions. Section 6 contains a numerical robustness study of the proposed approach. The adaptive scheme is compared with conventional adaptive controllers and its robustness with respect to vanishing and non-vanishing perturbations is studied. Finally, the concluding remarks are provided in Section 7.

All proofs of the main results are given in a companion paper [24]. The companion paper provides also a convenient algorithmic way of checking the parameter-observability assumption (H3) for a certain class of nonlinear control systems.

*Notation.*

* For a vector $x \in \Re^n$ we denote by $|x|$ its usual Euclidean norm, by $x'$ its transpose. For a real matrix $A \in \Re^{n \times m}$, $A' \in \Re^{m \times n}$ denotes its transpose and $|A| := \sup\{|Ax|; x \in \Re^n, |x| = 1\}$ is its induced norm. For a square matrix $A \in \Re^{n \times n}$, $\det(A)$ denotes its determinant.
* $\Re_+$ denotes the set of non-negative real numbers. $Z_+$ denotes the set of non-negative integers.
* We say that a function $V : \Re^n \to \Re_+$ is positive definite if $V(x) > 0$ for all $x \neq 0$ and $V(0) = 0$. We say that a continuous function $V : \Re^n \to \Re_+$ is radially unbounded if the following property holds: "for every $M > 0$ the set $\{x \in \Re^n : V(x) \leq M\}$ is compact".
* By $K$ we denote the class of strictly increasing continuous functions $a : \Re_+ \to \Re_+$ with $a(0) = 0$. By $K_\infty$ we denote the class of strictly increasing continuous functions $a : \Re_+ \to \Re_+$ with $a(0) = 0$ and $\lim_{s \to +\infty} a(s) = +\infty$. By $KL$ we denote the set of all continuous functions $\sigma = \sigma(s,t) : \Re_+ \times \Re_+ \to \Re_+$ with the properties: (i) for each $t \geq 0$ the mapping $\sigma(\cdot, t)$ is of class $K$; (ii) for each $s \geq 0$, the mapping $\sigma(s, \cdot)$ is non-increasing with $\lim_{t \to +\infty} \sigma(s,t) = 0$.

All stability notions used in this paper are the standard stability notions for time-invariant systems (see [25]).

## 2. Problem Formulation

Consider the system
$$\dot{x} = f(x,u) + g(x,u)\theta$$
$$x \in \Re^n, u \in \Re^m, \theta \in \Re^l \tag{2.1}$$
where $f : \Re^n \times \Re^m \to \Re^n$, $g : \Re^n \times \Re^m \to \Re^{n \times l}$ are smooth mappings with $f(0,0) = 0$, $g(0,0) = 0$ and $\theta \in \Re^l$ is a vector of constant but unknown parameters.

We suppose that there exist a smooth mapping $k : \Re^l \times \Re^n \to \Re^m$ with $k(\theta, 0) = 0$ for all $\theta \in \Re^l$, two families of continuous, positive definite and radially unbounded functions $V_\theta, Q_\theta : \Re^n \to \Re_+$ parameterized by $\theta \in \Re^l$ with the mappings $\Re^l \times \Re^n \ni (\theta, x) \to Q_\theta(x)$, $\Re^l \times \Re^n \ni (\theta, x) \to V_\theta(x)$ being continuous and such that the following assumptions hold.

**(H1)** *For each $\theta \in \Re^l$, $0 \in \Re^n$ is Globally Asymptotically Stable (GAS) for the closed-loop system*
$$\dot{x} = f(x, k(\theta, x)) + g(x, k(\theta, x))\theta \tag{2.2}$$
*Moreover, for every $\theta \in \Re^l$, $x_0 \in \Re^n$ the solution $x(t) \in \Re^n$ of (2.2) with initial condition $x(0) = x_0$ satisfies the inequality $V_\theta(x(t)) \leq Q_\theta(x_0)$ for all $t \geq 0$.*

**(H2)** *For every non-empty, compact set $\Theta \subset \Re^l$, the following property holds: "for every $M \geq 0$ there exists $R > 0$ such that the implication $V_\theta(x) \leq M, \theta \in \Theta \Rightarrow |x| \leq R$ holds".*

Assumption (H1) is a standard stabilizability assumption (necessary for all possible adaptive control design methodologies). For nonlinear systems, the design of a globally stabilizing state feedback law $u = k(\theta, x)$ is usually performed with the use of a Control Lyapunov Function (CLF, see [7,22,28,47] and references therein). In such cases, the CLF itself can be selected to be equal to both positive definite and radially unbounded functions $V_\theta, Q_\theta : \Re^n \to \Re_+$. In general, the Lyapunov function for system (2.2) (guaranteed to exist) can be selected to be equal to both positive definite and radially unbounded functions $V_\theta, Q_\theta : \Re^n \to \Re_+$. However, crude estimates may also be used for the selection of the functions $V_\theta, Q_\theta : \Re^n \to \Re_+$. For example, if for every $\theta \in \Re^l$ $0 \in \Re^n$ is Globally Exponentially Stable for (2.2), then there exists a family of constants $M_\theta, \omega_\theta > 0$ parameterized by $\theta \in \Re^l$, such that for every $\theta \in \Re^l$, $x_0 \in \Re^n$ the solution of the closed-loop system (2.2) with initial conditions $x(0) = x_0$ satisfies the estimate $|x(t)| \leq M_\theta \exp(-\omega_\theta t)|x_0|$ for all $t \geq 0$. In such cases, we may select $V_\theta(x) = |x|$ and $Q_\theta(x) = M_\theta|x|$, or $V_\theta(x) = |x|^2$ and $Q_\theta(x) = M_\theta^2 |x|^2$. Assumption (H2) is a technical assumption, which requires a "uniform" coercivity property for $V_\theta$ on compact sets of $\Re^l$.

In order to be able to estimate the vector of constant but unknown parameters $\theta \in \Re^l$, we need an additional technical assumption.

**(H3)** *There exists a positive integer $N$ such that the following implication holds:*

*"If there exist times $0 = \tau_0 < \tau_1 < ... < \tau_N$, vectors $\theta, d_0, ..., d_N \in \Re^l$ with $d_i \neq 0$ for $i = 0, ..., N$ and a right differentiable mapping $x \in C^0([0, \tau_N]; \Re^n) \cap C^1([0, \tau_N] \setminus \{\tau_0, ..., \tau_N\}; \Re^n)$ satisfying $\dot{x}(t) = f(x(t), k(\theta + d_i, x(t))) + g(x(t), k(\theta + d_i, x(t)))\theta$ for $t \in [\tau_i, \tau_{i+1})$, $i = 0, ..., N-1$, $g(x(t), k(\theta + d_j, x(t)))d_{i+1} = 0$ for all $t \in [\tau_j, \tau_{j+1}]$, $i = 0, ..., N-1$, $j = 0, ..., i$, then $x(t) = 0$ for all $t \in [0, \tau_N]$."*

Assumption (H3) is an observability assumption for the closed-loop system (2.1) with $u = k(\hat{\theta}, x)$, which guarantees that the only solution for which the vector of constant but unknown parameters $\theta \in \Re^l$ cannot be estimated is the zero solution. Assumption (H3) replaces the well-known "persistency of excitation" condition that is used in many cases for the design of adaptive control schemes. One of the advantages of assumption (H3) is that it can be verified a priori without additional assumptions: this feature is illustrated in the examples that follow.

## 3. Event-Triggered Identifier for a Certainty-Equivalence Adaptive Controller

In this section we gradually introduce the adaptive control law. The reader interested in a quick access to the adaptive controller may immediately refer to (3.1), (3.2), (3.4), (3.5), (3.10), (3.15) and then resume reading the rest of this section for explanations.

The control action in the interval between two consecutive events is governed by the nominal feedback $u = k(\theta, x)$ with the unknown $\theta \in \Re^l$ replaced by its estimate $\hat{\theta}$ at the beginning of the interval. Moreover, the estimate $\hat{\theta}$ of the unknown $\theta \in \Re^l$ is kept constant between two consecutive events. In other words, we have

$$u(t) = k(\hat{\theta}(\tau_i), x(t)) \quad , \quad t \in [\tau_i, \tau_{i+1}), i \in Z_+$$
$$\hat{\theta}(t) = \hat{\theta}(\tau_i) \quad , \quad t \in [\tau_i, \tau_{i+1}), i \in Z_+ \quad (3.1)$$

where $\{\tau_i \geq 0\}_{i=0}^{\infty}$ is the sequence of times of the events that satisfies

$$\tau_{i+1} = \min(\tau_i + T, r_i) \quad , \quad i \in Z_+$$
$$\tau_0 = 0 \quad (3.2)$$

where $T > 0$ is a positive constant (one of the tunable parameters of the proposed scheme) and $r_i > \tau_i$ is a time instant determined by the event trigger.

Let $a: \Re^n \to \Re_+$ be a continuous, positive definite function (again, one of the tunable parameters of the proposed scheme). The event trigger sets $r_i > \tau_i$ to be the smallest value of time $t > \tau_i$ for which

$$V_{\hat{\theta}(\tau_i)}(x(t)) = Q_{\hat{\theta}(\tau_i)}(x(\tau_i)) + a(x(\tau_i)) \qquad (3.3)$$

where $x(t)$ denotes the solution of (2.1) with $u(t) = k(\hat{\theta}(\tau_i), x(t))$. For the case that a time $t > \tau_i$ satisfying (3.3) does not exist, we set $r_i = +\infty$. For the case $x(\tau_i) = 0$ (and thus $a(x(\tau_i)) = 0$) we set $r_i := T$.

Formally, the event trigger is described by the equations:

$$r_i := \inf\left\{t > \tau_i : V_{\hat{\theta}(\tau_i)}(x(t)) = Q_{\hat{\theta}(\tau_i)}(x(\tau_i)) + a(x(\tau_i))\right\}, \text{ for } x(\tau_i) \neq 0 \qquad (3.4)$$

$$r_i := T, \text{ for } x(\tau_i) = 0 \qquad (3.5)$$

The description of the event-triggered adaptive control scheme is completed by the parameter update law, which will be activated at the times of the events.

In order to estimate the unknown vector $\theta \in \Re^l$, we notice that (by virtue of (2.1)) for every $t, \sigma \geq 0$ the following equation holds:

$$x(t) - x(\sigma) = \int_\sigma^t f(x(s), u(s))ds + \left(\int_\sigma^t g(x(s), u(s))ds\right)\theta \qquad (3.6)$$

Let $\tilde{N} > N$ be an (arbitrary; the last of the tunable parameters of the proposed scheme) positive integer that satisfies $\tilde{N} > N$, where $N > 0$ is the positive integer involved in Assumption (H3). Define for every $i \in Z_+$ the function $h_i : \Re^l \to \Re_+$ by the formula

$$h_i(\vartheta) := \int_{\mu_{i+1}}^{\tau_{i+1}} \int_{\mu_{i+1}}^{\tau_{i+1}} |p(t,\sigma) - q(t,\sigma)\vartheta|^2 d\sigma\, dt \qquad (3.7)$$

where

$$p(t,\sigma) := x(t) - x(\sigma) - \int_\sigma^t f(x(s), u(s))ds \qquad (3.8)$$

$$q(t,\sigma) := \int_\sigma^t g(x(s), u(s))ds \qquad (3.9)$$

$$\mu_{i+1} := \min\left\{\tau_j : j \in \{0,\ldots,i\}, \tau_j \geq \tau_{i+1} - \tilde{N}T\right\}. \qquad (3.10)$$

It follows from (3.6) and (3.7) that for every $i \in Z_+$ the function $h_i(\vartheta)$ has a global minimum at $\vartheta = \theta$ with $h_i(\theta) = 0$. Consequently, we get from Fermat's theorem that the following equation holds:

$$Z(\tau_{i+1}, \mu_{i+1}) = G(\tau_{i+1}, \mu_{i+1})\theta \qquad (3.11)$$

where

$$G(\tau_{i+1}, \mu_{i+1}) = \int_{\mu_{i+1}}^{\tau_{i+1}} \int_{\mu_{i+1}}^{\tau_{i+1}} q'(t,\sigma) q(t,\sigma) d\sigma\, dt$$

$$Z(\tau_{i+1}, \mu_{i+1}) = \int_{\mu_{i+1}}^{\tau_{i+1}} \int_{\mu_{i+1}}^{\tau_{i+1}} q'(t,\sigma) p(t,\sigma) d\sigma\, dt \qquad (3.12)$$

It should be noticed that the matrix $G(\tau_{i+1}, \mu_{i+1}) \in \Re^{l \times l}$ is symmetric and positive semi-definite. Consequently, if $G(\tau_{i+1}, \mu_{i+1}) \in \Re^{l \times l}$ is invertible (i.e., $\det(G(\tau_{i+1}, \mu_{i+1})) \neq 0$) then $G(\tau_{i+1}, \mu_{i+1}) \in \Re^{l \times l}$ is positive definite with $\det(G(\tau_{i+1}, \mu_{i+1})) > 0$ and

$$\theta = (G(\tau_{i+1}, \mu_{i+1}))^{-1} Z(\tau_{i+1}, \mu_{i+1}) \qquad (3.13)$$

Estimate (3.13) is nothing else but the least squares estimate of the unknown vector $\theta \in \Re^l$ on the interval $[\mu_{i+1}, \tau_{i+1}]$ (and is a modification of the dead-beat observer proposed in [21]). In the general case, the following quadratic optimization problem with linear equality constraints

$$\min\left\{\left|\vartheta - \hat{\theta}(\tau_i)\right|^2 : \vartheta \in \Re^l, Z(\tau_{i+1}, \mu_{i+1}) = G(\tau_{i+1}, \mu_{i+1})\vartheta\right\} \quad (3.14)$$

has a unique solution (which may or may not be equal to $\theta$ depending on whether $\det(G(\tau_{i+1}, \mu_{i+1})) > 0$ or $\det(G(\tau_{i+1}, \mu_{i+1})) = 0$).

We can therefore define the following parameter update law:

$$\hat{\theta}(\tau_{i+1}) = \arg\min\left\{\left|\vartheta - \hat{\theta}(\tau_i)\right|^2 : \vartheta \in \Re^l, Z(\tau_{i+1}, \mu_{i+1}) = G(\tau_{i+1}, \mu_{i+1})\vartheta\right\} \quad (3.15)$$

Equation (3.15) implies that $\hat{\theta}(\tau_{i+1})$ is the projection of $\hat{\theta}(\tau_i)$ on the hyperplane in the parameter space ($\Re^l$) defined by the linear equality constraints $Z(\tau_{i+1}, \mu_{i+1}) = G(\tau_{i+1}, \mu_{i+1})\vartheta$. It should be emphasized that $\hat{\theta}(\tau_{i+1})$ as given in (3.15) may or may not be dependent on $\hat{\theta}(\tau_i)$ (depending on the rank of the matrix $G(\tau_{i+1}, \mu_{i+1}) \in \Re^{l \times l}$). Moreover, it should also be noticed that the operator involved in (3.15) is not a continuous operator. However, in practice an accurate continuous approximation of the parameter update law (3.15) may be used (for example, by using an appropriate Tikhonov regularization procedure as the replacement of the linear equality constraints $Z(\tau_{i+1}, \mu_{i+1}) = G(\tau_{i+1}, \mu_{i+1})\vartheta$ by $Z(\tau_{i+1}, \mu_{i+1}) = (\eta I + G(\tau_{i+1}, \mu_{i+1}))\vartheta$, where $I \in \Re^{l \times l}$ denotes the unit matrix and $\eta > 0$ is a small constant; see [50]).

**Remark 3.1:** It is important to notice that the parameter update law given by (3.15) can be implemented by a set of additional ODEs. Indeed, an implementation of the parameter update law (3.15) is given by the following $(2n+l)(1+l)$ ODEs

$$\begin{aligned}
\dot{z} &= f(x,u) \quad , \quad z \in \Re^n \\
\dot{w} &= x - z \quad , \quad w \in \Re^n \\
\dot{B} &= g(x,u) \quad , \quad B \in \Re^{n \times l} \\
\dot{\varphi} &= B'(x-z) \quad , \quad \varphi \in \Re^l \\
\dot{Q} &= B \quad , \quad Q \in \Re^{n \times l} \\
\dot{R} &= B'B \quad , \quad R \in \Re^{l \times l}
\end{aligned} \quad (3.16)$$

with initial conditions $z(0) = w(0) = 0$, $B(0) = Q(0) = 0$, $\varphi(0) = 0$, $R(0) = 0$. The parameter update law (3.15) is given by (3.15), where

$$Z(\tau_{i+1}, \mu_{i+1}) := (\tau_{i+1} - \mu_{i+1})(\varphi(\tau_{i+1}) - \varphi(\mu_{i+1})) - (Q(\tau_{i+1}) - Q(\mu_{i+1}))'(w(\tau_{i+1}) - w(\mu_{i+1}))$$

$$G(\tau_{i+1}, \mu_{i+1}) := (\tau_{i+1} - \mu_{i+1})(R(\tau_{i+1}) - R(\mu_{i+1})) - (Q(\tau_{i+1}) - Q(\mu_{i+1}))'(Q(\tau_{i+1}) - Q(\mu_{i+1}))$$

However, it should be emphasized that in many cases, the structure of the control system (2.1) allows a large reduction of the number of ODEs that are needed for the implementation of the parameter update law given by (3.15).

**Remark 3.2:** A simpler least-squares identifier than (3.15) can be obtained by using the function $\tilde{h}_i(\vartheta) := \int_{\mu_{i+1}}^{\tau_{i+1}} |p(t, \mu_{i+1}) - q(t, \mu_{i+1})\vartheta|^2 dt$. In this way, by repeating the same arguments, we end up with the following parameter update law:

$$\hat{\theta}(\tau_{i+1}) = \arg\min\left\{\left|\vartheta - \hat{\theta}(\tau_i)\right|^2 : \vartheta \in \Re^l, \tilde{Z}(\tau_{i+1}, \mu_{i+1}) = \tilde{G}(\tau_{i+1}, \mu_{i+1})\vartheta\right\} \quad (3.17)$$

where

$$\tilde{G}(\tau_{i+1}, \mu_{i+1}) = \int_{\mu_{i+1}}^{\tau_{i+1}} q'(t, \mu_{i+1}) q(t, \mu_{i+1}) dt$$

$$\tilde{Z}(\tau_{i+1}, \mu_{i+1}) = \int_{\mu_{i+1}}^{\tau_{i+1}} q'(t, \mu_{i+1}) p(t, \mu_{i+1}) dt$$

This is exactly the parameter update law that would be derived by using the dead-beat estimator in [21] without modification. The parameter update law (3.17) is simpler than the update (3.15), since it involves simple integrals instead of double integrals. It can be used without any problem and gives the same results that are presented in the following section. However, there is a big difference between (3.17) and (3.15): (3.17) gives much more weight to the measurement $x(\mu_{i+1})$ compared to the weight given to measurements in other times, while (3.15) gives equal weight to all measurements. It is therefore expected that (3.17) will work better than (3.15) only in cases where the measurements at the times of the events are to be trusted more than measurements in times between events. In all other cases, where the measurements can be trusted "equally", (3.15) is expected to be more robust than (3.17) with respect to random measurement noise.

## 4. Statements of Stability Results

We consider the plant (2.1) with the controller (3.1), (3.2), (3.4), (3.5) and the parameter estimator (3.15), (3.10). The first main result guarantees global regulation of $x$ to zero and a dead-beat estimation of the unknown vector $\theta \in \Re^l$ for all non-zero initial conditions.

**Theorem 4.1:** *Consider the control system (2.1) under assumptions (H1), (H2), (H3). Let $T > 0$ be a positive constant and let $a: \Re^n \to \Re_+$ be a continuous, positive definite function. Finally, let $\tilde{N} > N$ be a positive integer that satisfies $\tilde{N} > N$, where $N > 0$ is the positive integer involved in Assumption (H3). Then there exists a family of KL mappings $\tilde{\sigma}_{\theta, \hat{\theta}} \in KL$ parameterized by $\theta \in \Re^l$, $\hat{\theta} \in \Re^l$ such that for every $\theta \in \Re^l$, $x_0 \in \Re^n$, $\hat{\theta}_0 \in \Re^l$ the solution of the hybrid closed-loop system (2.1) with (3.1), (3.2), (3.4), (3.5), (3.10), (3.15) and initial conditions $x(0) = x_0$, $\hat{\theta}(0) = \hat{\theta}_0$ is unique, is defined for all $t \geq 0$ and satisfies $|x(t)| \leq \tilde{\sigma}_{\theta, \hat{\theta}_0}(|x_0|, t)$ for all $t \geq 0$. Moreover, if $x_0 \neq 0$ then $\hat{\theta}(t) = \theta$ for all $t \geq NT$.*

The second main result guarantees local exponential regulation of $x$ to zero under the assumption that the nominal feedback law $u = k(\theta, x)$ achieves local exponential stabilization.

**Theorem 4.2:** *Consider the control system (2.1) under assumptions (H1), (H2), (H3). Moreover, suppose that for each $\theta \in \Re^l$, there exist constants $M_\theta, \omega_\theta, R_\theta > 0$ such that for every $x_0 \in \Re^n$ with $|x_0| \leq R_\theta$ the solution of (2.2) with initial condition $x(0) = x_0$ satisfies the estimate $|x(t)| \leq M_\theta \exp(-\omega_\theta t)|x_0|$ for all $t \geq 0$; i.e., $0 \in \Re^n$ is Locally Exponentially Stable (LES) for the closed-loop system (2.2). Furthermore, suppose that for every nonempty, compact set $\Theta \subset \Re^l$ there exist constants $R > 0$, $K_2 > K_1 > 0$ such that*

$$K_1 |x|^2 \leq V_\theta(x) \leq Q_\theta(x) \leq K_2 |x|^2,$$
*for all $x \in \Re^n, \theta \in \Theta$ with $|x| \leq R$* (4.1)

*Let $T > 0$ be a positive constant and let $a: \Re^n \to \Re_+$ be a continuous, positive definite function that satisfies $\sup\{|x|^{-2} a(x) : x \in \Re^n, x \neq 0, |x| \leq \delta\} < +\infty$ for certain $\delta > 0$. Finally, let $\tilde{N} > N$ be a positive integer that satisfies $\tilde{N} > N$, where $N > 0$ is the positive integer involved in Assumption (H3). Then there exists a family of constants $\tilde{M}_{\theta, \hat{\theta}}, \tilde{R}_{\theta, \hat{\theta}} > 0$ parameterized by $(\theta, \hat{\theta}) \in \Re^l \times \Re^l$, such that for every $\theta \in \Re^l$, $x_0 \in \Re^n$, $\hat{\theta}_0 \in \Re^l$ with $|x_0| \leq \tilde{R}_{\theta, \hat{\theta}_0}$ the solution of the hybrid closed-loop system (2.1) with* (3.1), (3.2), (3.4), (3.5), (3.10), (3.15) *and initial conditions* $x(0) = x_0$, $\hat{\theta}(0) = \hat{\theta}_0$ *satisfies the estimate* $|x(t)| \leq \tilde{M}_{\theta, \hat{\theta}_0} \exp(-\omega_\theta t)|x_0|$ *for all* $t \geq 0$.

It should be noticed that Theorem 4.2 guarantees that the local exponential stability estimate $|x(t)| \leq \tilde{M}_{\theta, \hat{\theta}_0} \exp(-\omega_\theta t)|x_0|$ holds when $|x_0| \leq \tilde{R}_{\theta, \hat{\theta}_0}$ and for arbitrary initial condition $\hat{\theta}_0 \in \Re^l$. In other words, the adjective "local" refers only to $x$ and not to $\hat{\theta}$. Moreover, the reader should notice that the event-triggered adaptive scheme (3.1), (3.2), (3.4), (3.5), (3.10), (3.15) guarantees convergence with the same convergence rate $\omega_\theta$ as the nominal feedback controller with known parameter values. Finally, it should be emphasized that in addition to the exponential stability estimate $|x(t)| \leq \tilde{M}_{\theta, \hat{\theta}_0} \exp(-\omega_\theta t)|x_0|$, Theorem 4.2 guarantees all the conclusions of Theorem 4.1 (because all assumptions of Theorem 4.1 are fulfilled).

The following result guarantees global exponential regulation of $x$ to zero under the assumption that the nominal feedback law $u = k(\theta, x)$ achieves global exponential stabilization.

**Theorem 4.3:** *Consider system (2.1) under assumptions (H1), (H2), (H3). Moreover, suppose that for each* $\theta \in \Re^l$, $0 \in \Re^n$ *is Globally Exponentially Stable (GES) for the closed-loop system (2.2) and that for every nonempty, compact set* $\Theta \subset \Re^l$ *there exist constants* $K_2 > K_1 > 0$ *such that*

$$K_1|x|^2 \leq V_\theta(x) \leq Q_\theta(x) \leq K_2|x|^2,$$
*for all* $x \in \Re^n, \theta \in \Theta$ \hfill (4.2)

*Let* $T > 0$ *be a positive constant and let* $a: \Re^n \to \Re_+$ *be a continuous, positive definite function that satisfies* $\sup\{|x|^{-2} a(x) : x \in \Re^n, x \neq 0\} < +\infty$. *Finally, let* $\tilde{N} > N$ *be a positive integer that satisfies* $\tilde{N} > N$, *where* $N > 0$ *is the positive integer involved in Assumption (H3). Then there exists a family of constants* $\tilde{M}_{\theta, \hat{\theta}} > 0$ *parameterized by* $\theta \in \Re^l$, $\hat{\theta} \in \Re^l$, *such that for every* $\theta \in \Re^l$, $x_0 \in \Re^n$, $\hat{\theta}_0 \in \Re^l$ *the solution of the hybrid closed-loop system (2.1) with (3.1), (3.2), (3.4), (3.5), (3.10), (3.15) and initial conditions* $x(0) = x_0$, $\hat{\theta}(0) = \hat{\theta}_0$ *satisfies the estimate* $|x(t)| \leq \tilde{M}_{\theta, \hat{\theta}_0} \exp(-\omega_\theta t)|x_0|$ *for all* $t \geq 0$.

Finally, the following corollary deals with the case of Linear Time-Invariant (LTI) systems with unknown parameters.

**Corollary 4.4:** *Consider the system*

$$\dot{x} = (A + \theta_1 C_1 + ... + \theta_l C_l)x + Bu$$
$$x \in \Re^n, u \in \Re^m, \theta = (\theta_1, ..., \theta_l)' \in \Re^l$$ \hfill (4.3)

*where* $A, C_1, ..., C_l \in \Re^{n \times n}$, $B \in \Re^{n \times m}$ *are constant matrices. Suppose that there exists a family of constants* $\omega_\theta > 0$ *parameterized by* $\theta \in \Re^l$ *and a continuous mapping* $\Re^l \ni \theta \to M(\theta) \in [1, +\infty)$ *such that* $|\exp(t(A + \theta_1 C_1 + ... + \theta_l C_l + BK_\theta))| \leq \exp(-\omega_\theta t)M(\theta)$ *for all* $t \geq 0$. *Moreover, suppose that for every* $\theta = (\theta_1, ..., \theta_l)' \in \Re^l$, $\hat{\theta} \in \Re^l$, $\vartheta = (\vartheta_1, ..., \vartheta_l)' \in \Re^l$ *with* $\hat{\theta} \neq \theta$ *and* $\vartheta \neq 0$, *the pair of matrices* $(A + \theta_1 C_1 + ... + \theta_l C_l + BK_{\hat{\theta}}, \vartheta_1 C_1 + ... + \vartheta_l C_l)$ *is an observable pair of matrices. Let* $a, T > 0$ *be constants and let* $\tilde{N} > 1$ *be a positive integer. Let* $L: \Re^n \to \Re^{n \times l}$ *be the linear operator defined by* $L * x = (C_1 x e_1' + ... + C_l x e_l') \in \Re^{n \times l}$ *for* $x \in \Re^n$ *with* $e_1' = (1, 0, ..., 0)' \in \Re^l$, ..., $e_l' = (0, ..., 0, 1)' \in \Re^l$. *Then there exists a family of constants* $\tilde{M}_{\theta, \hat{\theta}} > 0$ *parameterized by* $\theta \in \Re^l$, $\hat{\theta} \in \Re^l$, *such that for every* $\theta \in \Re^l$, $x_0 \in \Re^n$, $\hat{\theta}_0 \in \Re^l$ *the solution of the hybrid closed-loop system (4.3) with (3.2), (3.5), (3.10),*

$$u(t) = K_{\hat{\theta}(\tau_i)} x(t) \quad , \quad t \in [\tau_i, \tau_{i+1}), i \in Z_+$$
$$\hat{\theta}(t) = \hat{\theta}(\tau_i) \quad , \quad t \in [\tau_i, \tau_{i+1}), i \in Z_+$$ \hfill (4.4)

$$r_i := \inf\left\{t > \tau_i : |x(t)| = |x(\tau_i)|\sqrt{a + M^2(\hat{\theta}(\tau_i))}\right\}, \text{ for } x(\tau_i) \neq 0, \tag{4.5}$$

$$\hat{\theta}(\tau_{i+1}) = \arg\min\left\{\left|\vartheta - \hat{\theta}(\tau_i)\right|^2 : \vartheta \in \Re^l, q(\tau_{i+1}, \mu_{i+1}) = Q(\tau_{i+1}, \mu_{i+1})\vartheta\right\}, \tag{4.6}$$

*where*

$$\begin{aligned}
&\dot{z} = x, \; z \in \Re^n \\
&\dot{w} = u, \; w \in \Re^m \\
&y = x - Az - Bw \\
&q(\tau, \mu) = \int_\mu^\tau \int_\mu^\tau (L*(z(t) - z(\sigma)))'(y(t) - y(\sigma))\,d\sigma\,dt \\
&Q(\tau, \mu) = \int_\mu^\tau \int_\mu^\tau (L*(z(t) - z(\sigma)))'(L*(z(t) - z(\sigma)))\,d\sigma\,dt
\end{aligned} \tag{4.7}$$

*and initial conditions* $x(0) = x_0$, $\hat{\theta}(0) = \hat{\theta}_0$, $z(0) = 0$, $w(0) = 0$ *satisfies the estimate* $|x(t)| \leq \tilde{M}_{\theta,\hat{\theta}_0} \exp(-\omega_\theta t)|x_0|$ *for all* $t \geq 0$. *Moreover, if* $x_0 \neq 0$ *then* $\hat{\theta}(t) = \theta$ *for all* $t \geq T$.

## 5. Example: Illustration of the Assumptions and Design

As remarked above, assumptions (H1), (H2), (H3) can be checked a priori without knowledge of the solutions of system (2.1). This feature is illustrated by the following example, which deals with a system not of the form (2.3).

Consider the following nonlinear, planar system

$$\begin{aligned}
\dot{x}_1 &= \theta_1 x_1 + \theta_2 x_1^2 + x_2 \\
\dot{x}_2 &= u \\
x &= (x_1, x_2)' \in \Re^2, u \in \Re, \theta = (\theta_1, \theta_2)' \in \Re^2
\end{aligned} \tag{5.1}$$

System (5.1) is a system of the form (2.1) with

$$f(x,u) = \begin{bmatrix} x_2 \\ u \end{bmatrix}, \; g(x,u) = \begin{bmatrix} x_1 & x_1^2 \\ 0 & 0 \end{bmatrix}, \text{ for all } (x,u) \in \Re^2 \times \Re. \tag{5.2}$$

Using backstepping we are in a position to design a stabilizing feedback for system (3.3). More specifically, the smooth feedback law:

$$\begin{aligned}
k(\theta, x) &:= -x_1 - (1 + \theta_1 + 2\theta_2 x_1)(x_2 + \theta_1 x_1 + \theta_2 x_1^2) \\
&\quad - (x_2 + x_1 + \theta_1 x_1 + \theta_2 x_1^2) \\
&\text{for all } (x, \theta) \in \Re^2 \times \Re^2
\end{aligned} \tag{5.3}$$

guarantees that the derivative of the CLF

$$V_\theta(x) := \frac{1}{2}x_1^2 + \frac{1}{2}(x_2 + x_1 + \theta_1 x_1 + \theta_2 x_1^2)^2, \text{ for all } (x, \theta) \in \Re^2 \times \Re^2 \tag{5.4}$$

along the solutions of the closed-loop system (5.1) with $u = k(\theta, x)$ satisfies

$$\dot{V}_\theta(x) = -x_1^2 - (x_2 + x_1 + \theta_1 x_1 + \theta_2 x_1^2)^2 = -2V_\theta(x), \text{ for all } (x,\theta) \in \Re^2 \times \Re^2 \tag{5.5}$$

Therefore, Assumption (H1) holds with $Q_\theta \equiv V_\theta$.

In order to show the validity of Assumption (H2), we consider an arbitrary non-empty, compact set $\Theta \subset \Re^2$. Let $\rho > 0$ be sufficiently large so that the implication $\theta \in \Theta \Rightarrow |\theta| \leq \rho$ holds. Let $M \geq 0$ be given and suppose that $V_\theta(x) \leq M, \theta \in \Theta$. Definition (5.4) implies that $|x_2 + x_1 + \theta_1 x_1 + \theta_2 x_1^2| \leq \sqrt{2M}$, $|x_1| \leq \sqrt{2M}$. Using the triangle inequality $|x_2| \leq |x_2 + x_1 + \theta_1 x_1 + \theta_2 x_1^2| + |x_1 + \theta_1 x_1 + \theta_2 x_1^2|$ and the fact that $|\theta| \leq \rho$, we obtain by repeated use of the triangle inequality:

$$|x_2| \leq (2+\rho)\sqrt{2M} + 2\rho M$$

Consequently, we get $|x| \leq |x_1| + |x_2| \leq R = (3+\rho)\sqrt{2M} + 2\rho M$ and therefore Assumption (H2) holds.

Next we show that Assumption (H3) holds with $N = 2$. Let times $0 = \tau_0 < \tau_1 < \tau_2$, vectors $\theta, d_0, d_1, d_2 \in \Re^2$ with $d_i \neq 0$ for $i = 0,1,2$ and a right differentiable mapping $x \in C^0([0, \tau_2]; \Re^2) \cap C^1((0, \tau_2) \setminus \{\tau_1\}; \Re^2)$ satisfying

$$\dot{x}(t) = f(x(t), k(\theta + d_0, x(t))) + g(x(t), k(\theta + d_0, x(t)))\theta, \text{ for } t \in [\tau_0, \tau_1), \quad (5.6)$$

$$\dot{x}(t) = f(x(t), k(\theta + d_1, x(t))) + g(x(t), k(\theta + d_1, x(t)))\theta, \text{ for } t \in [\tau_1, \tau_2), \quad (5.7)$$

$$g(x(t), k(\theta + d_0, x(t)))d_1 = 0, \text{ for all } t \in [0, \tau_1], \quad (5.8)$$

$$g(x(t), k(\theta + d_0, x(t)))d_2 = 0, \text{ for all } t \in [0, \tau_1], \quad (5.9)$$

$$g(x(t), k(\theta + d_1, x(t)))d_2 = 0, \text{ for all } t \in [\tau_1, \tau_2], \quad (5.10)$$

be given. Equalities (5.2), (5.8), (5.9), (5.10) imply that

$$d_{1,1} x_1(t) + d_{1,2} x_1^2(t) = 0, \text{ for all } t \in [0, \tau_1] \quad (5.11)$$

$$d_{2,1} x_1(t) + d_{2,2} x_1^2(t) = 0, \text{ for all } t \in [0, \tau_2] \quad (5.12)$$

where $d_1 = (d_{1,1}, d_{1,2})' \neq 0$, $d_2 = (d_{2,1}, d_{2,2})' \neq 0$. Equalities (5.11), (5.12) imply the existence of a constant $c \in \Re$ with $d_{1,1}c + d_{1,2}c^2 = d_{2,1}c + d_{2,2}c^2 = 0$ such that $x_1(t) = c$ for all $t \in [0, \tau_2]$. Consequently, we must have $\dot{x}_1(t) = 0$ for all $t \in [0, \tau_2]$, which implies that $x_2(t) = -\theta_1 c - \theta_2 c^2$ for all $t \in [0, \tau_2]$. Therefore, we must have $\dot{x}_2(t) = 0$ for all $t \in [0, \tau_2]$, which, combined with (5.3), (5.7) and the facts $x_1(t) = c$, $x_2(t) = -\theta_1 c - \theta_2 c^2$, $d_{1,1}c + d_{1,2}c^2 = 0$, implies that $c = 0$. It follows that $x(t) = 0$ for all $t \in [0, \tau_2]$.

It should be also remarked that for each $\theta \in \Re^2$, $0 \in \Re^2$ is LES for the closed-loop system (5.1) with $u = k(\theta, x)$. Furthermore, for every nonempty, compact set $\Theta \subset \Re^2$ there exist constants $R > 0$, $K_2 > K_1 > 0$ such that (4.1) holds. Indeed, let an arbitrary non-empty, compact set $\Theta \subset \Re^2$ be given. Let $\rho > 0$ be sufficiently large so that the implication $\theta \in \Theta \Rightarrow |\theta| \leq \rho$ holds. Using the inequalities $(1-\varepsilon)a^2 + (1-\varepsilon^{-1})b^2 \leq (a+b)^2 \leq 2a^2 + 2b^2$ (that hold for all $a, b \in \Re$ and $\varepsilon > 0$) and definition (5.4) we get for all $(\theta, x) \in \Re^2 \times \Re^2$:

$$\frac{1}{2}\left(1 + (1-\varepsilon)(1 + \theta_1 + \theta_2 x_1)^2\right)x_1^2 + \frac{1}{2}(1-\varepsilon^{-1})x_2^2 \leq V_\theta(x) \leq x_2^2 + \frac{1}{2}x_1^2 + x_1^2(1 + \theta_1 + \theta_2 x_1)^2$$

Let $R > 0$ be given. The above inequalities with $\varepsilon = 1 + \dfrac{1}{2(1 + \rho + \rho R)^2}$ in conjunction with the implication $\theta \in \Theta \Rightarrow |\theta| \leq \rho$, imply that for all $(\theta, x) \in \Theta \times \Re^2$ with $|x| \leq R$, it holds that

$$\frac{1}{4}x_1^2 + \frac{1}{4(1 + \rho + \rho R)^2 + 2}x_2^2 \leq V_\theta(x) \leq x_2^2 + \left(\frac{1}{2} + (1 + \rho + \rho R)^2\right)x_1^2$$

Using the above inequality and the fact that $Q_\theta \equiv V_\theta$, we conclude that (4.1) holds with $K_1 = \dfrac{1}{4(1 + \rho + \rho R)^2 + 2}$ and $K_2 = (1 + \rho + \rho R)^2 + 1/2$.

It follows that all assumptions of Theorem 4.1 and Theorem 4.2 hold. Therefore, for every $T > 0$, for every continuous, positive definite function $a: \Re^n \to \Re_+$ and for every $\tilde{N} > 2$ there exists a family of $KL$ mappings $\tilde{\sigma}_{\theta,\hat{\theta}} \in KL$ parameterized by $\theta \in \Re^2$, $\hat{\theta} \in \Re^2$ such that for every $\theta \in \Re^2$, $x_0 \in \Re^2$, $\hat{\theta}_0 \in \Re^2$ the solution of the hybrid closed-loop system (5.1) with (3.1), (3.2), (3.4), (3.5),

$$\hat{\theta}(\tau_{i+1}) = \arg\min\left\{ \left|\vartheta - \hat{\theta}(\tau_i)\right|^2 : A = C_1 \vartheta_1 + C_2 \vartheta_2, B = C_2 \vartheta_1 + C_3 \vartheta_2 \right\} \tag{5.13}$$

with

$$A = \int_{\mu_{i+1}}^{\tau_{i+1}} \int_{\mu_{i+1}}^{\tau_{i+1}} \left( x_1(t) - x_1(\sigma) - \int_\sigma^t x_2(s)ds \right) \left( \int_\sigma^t x_1(s)ds \right) d\sigma\, dt \quad,\quad B = \int_{\mu_{i+1}}^{\tau_{i+1}} \int_{\mu_{i+1}}^{\tau_{i+1}} \left( x_1(t) - x_1(\sigma) - \int_\sigma^t x_2(s)ds \right) \left( \int_\sigma^t x_1^2(s)ds \right) d\sigma\, dt$$

$$C_1 = \int_{\mu_{i+1}}^{\tau_{i+1}} \int_{\mu_{i+1}}^{\tau_{i+1}} \left( \int_\sigma^t x_1(s)ds \right)^2 d\sigma\, dt \quad,\quad C_2 = \int_{\mu_{i+1}}^{\tau_{i+1}} \int_{\mu_{i+1}}^{\tau_{i+1}} \left( \int_\sigma^t x_1(s)ds \right)\left( \int_\sigma^t x_1^2(s)ds \right) d\sigma\, dt \quad,\quad C_3 = \int_{\mu_{i+1}}^{\tau_{i+1}} \int_{\mu_{i+1}}^{\tau_{i+1}} \left( \int_\sigma^t x_1^2(s)ds \right)^2 d\sigma\, dt$$

$$\mu_{i+1} := \min\left\{ \tau_j : j \in \{0,\ldots,i\}, \tau_j \geq \tau_{i+1} - \tilde{N}T \right\}$$

initial conditions $x(0) = x_0$, $\hat{\theta}(0) = \hat{\theta}_0$ is unique, is defined for all $t \geq 0$ and satisfies $|x(t)| \leq \tilde{\sigma}_{\theta,\hat{\theta}_0}(|x_0|, t)$ for all $t \geq 0$. Moreover, if $x_0 \neq 0$ then $\hat{\theta}(t) = \theta$ for all $t \geq 2T$. Finally, if there exists $\delta > 0$ with $\sup\left\{ |x|^{-2} a(x) : x \neq 0, |x| \leq \delta \right\} < +\infty$, then there exists a family of constants $\tilde{M}_{\theta,\hat{\theta}}, \omega_\theta, \tilde{R}_{\theta,\hat{\theta}} > 0$ parameterized by $\theta \in \Re^2$, $\hat{\theta} \in \Re^2$, such that for every $\theta \in \Re^2$, $x_0 \in \Re^2$, $\hat{\theta}_0 \in \Re^2$ with $|x_0| \leq \tilde{R}_{\theta,\hat{\theta}_0}$ the solution of the hybrid closed-loop system (5.1) with (3.1), (3.2), (3.4), (3.5), (5.13) and initial conditions $x(0) = x_0$, $\hat{\theta}(0) = \hat{\theta}_0$ satisfies the estimate $|x(t)| \leq \tilde{M}_{\theta,\hat{\theta}_0} \exp(-\omega_\theta t) |x_0|$ for all $t \geq 0$.

## 6. Robustness Tests

Consider the system

$$\begin{aligned} \dot{x}_1 &= (\theta + v_1) x_1^2 + x_2 + v_2 \\ \dot{x}_2 &= u \\ x &= (x_1, x_2) \in \Re^2, u \in \Re, v = (v_1, v_2) \in \Re^2 \end{aligned} \tag{6.1}$$

where $x \in \Re^2$ is the state, $\theta \in \Re$ is an unknown constant parameter, $u \in \Re$ is the control input and $v = (v_1, v_2) \in \Re^2$ is an unknown, time-varying disturbance. Our goal is to test the robustness properties of the closed-loop system (6.1) with respect to the disturbances $v = (v_1, v_2) \in \Re^2$ when the proposed event-triggered adaptive control scheme is used. The performance of the event-triggered adaptive control scheme is compared with conventional adaptive schemes and the nominal feedback law with known parameter $\theta \in \Re$.

Using backstepping we are in a position to design a stabilizing feedback for system (6.1). More specifically, the smooth feedback law:

$$\begin{aligned} k(\theta, x) &:= -x_1 - \left(1 + 2\theta x_1 + 3x_1^2\right)\left(\theta x_1^2 + x_2\right) \\ &\quad - \frac{1}{2}\left( x_2 + x_1 + x_1^3 + \theta x_1^2 \right)\left( 1 + \left(1 + 2\theta x_1 + 3x_1^2\right)^2 \left(1 + x_1^4\right) \right) \end{aligned} \tag{6.2}$$

guarantees:
i) Global Asymptotic Stability and Local Exponential Stability for the disturbance-free closed-loop system (6.1) with $u = k(\theta, x)$,
ii) Robust Global Asymptotic Stability w.r.t. $v_1$ for the closed-loop system (6.1) with $u = k(\theta, x)$ and $v_2 \equiv 0$, provided that $|v_1| \leq c$, where $c > 0$ is a sufficiently small constant,
iii) Input-to-State Stability w.r.t. $v \in \Re^2$ for the closed-loop system (6.1) with $u = k(\theta, x)$.
More specifically, the above properties are proved by using the Lyapunov function

$$V_\theta(x) := \frac{1}{2} x_1^2 + \frac{1}{2}\left( x_2 + x_1 + x_1^3 + \theta x_1^2 \right)^2 \tag{6.3}$$

which satisfies the inequality $\dot{V}_\theta(x) \leq -V_\theta(x) + \frac{2+v_1^2}{4}v_1^2 + v_2^2$ for all $x \in \mathcal{R}^2$, $v \in \mathcal{R}^2$, along the solutions of the closed-loop system (6.1) with $u = k(\theta, x)$.

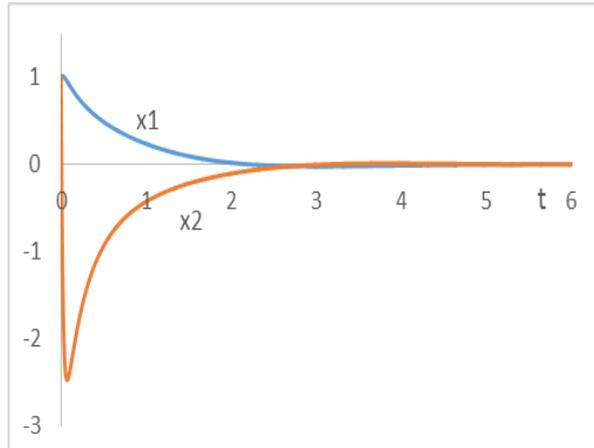

**Fig. 1:** The closed-loop system (6.1) with $u = k(\theta, x)$, $A_1 = A_2 = 0$.

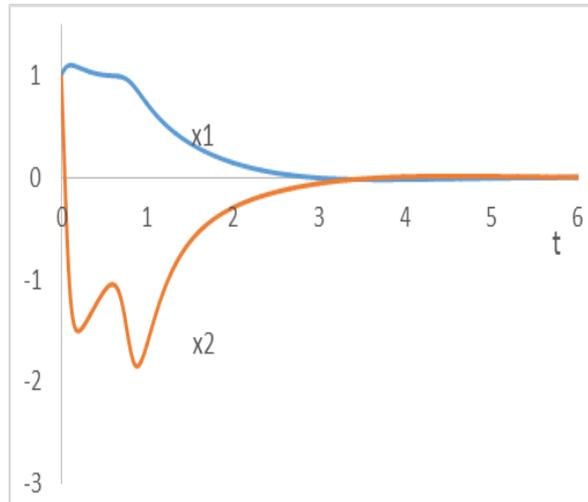

**Fig. 2:** The closed-loop system (6.1) with (6.4), $A_1 = A_2 = 0$.

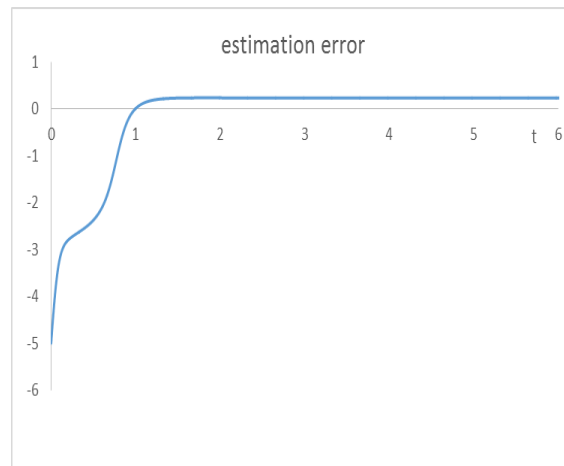

**Fig. 3:** The estimation error $\hat{\theta}(t) - \theta$ for the closed-loop system (6.1) with (6.4), $A_1 = A_2 = 0$.

Next, we design an adaptive control law for the disturbance-free system (6.1). We are using the extended matching approach (see pages 124-127 of the book [28]), the Lyapunov function given by (6.3) and the nominal feedback law given by (6.2). In this way, we get the dynamic regulator:

$$\frac{d\hat{\theta}}{dt} = \gamma x_1^2 \left( x_1 + \left( x_2 + x_1 + x_1^3 + \hat{\theta} x_1^2 \right)\left( 1 + 2\hat{\theta} x_1 + 3x_1^2 \right) \right)$$

$$u = k(\hat{\theta}, x) - \gamma x_1^4 \left( x_1 + \left( x_2 + x_1 + x_1^3 + \hat{\theta} x_1^2 \right)\left( 1 + 2\hat{\theta} x_1 + 3x_1^2 \right) \right)$$

(6.4)

where $\gamma > 0$ is a constant.

Finally, we design an event-triggered adaptive control scheme for the disturbance-free system (6.1) based on the Lyapunov function given by (6.3) and the nominal feedback law given by (6.2). Indeed, following the same procedure as the one followed in the example of Section 5, we are in a position to show that Assumptions (H1), (H2), (H3) hold with $Q_\theta \equiv V_\theta$ and $N = 1$. Furthermore, for every nonempty, compact set $\Theta \subset \Re$ there exist constants $R > 0$, $K_2 > K_1 > 0$ such that (4.1) holds. It follows that all assumptions of Theorem 4.1 and Theorem 4.2 hold.

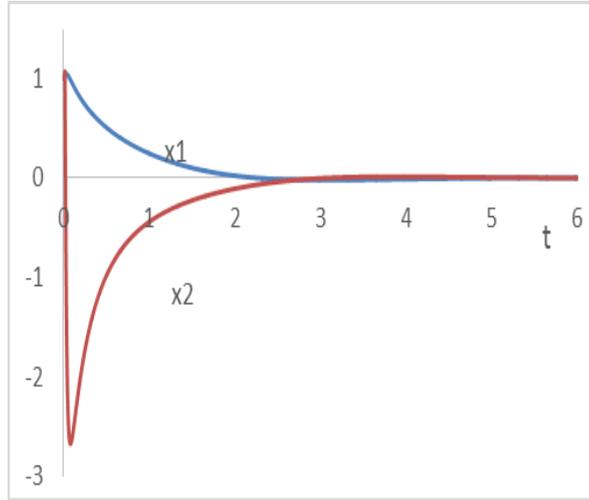

**Fig. 4:** The closed-loop system (6.1) with (6.5), (6.6), (6.12), (6.13), $A_1 = A_2 = 0$.

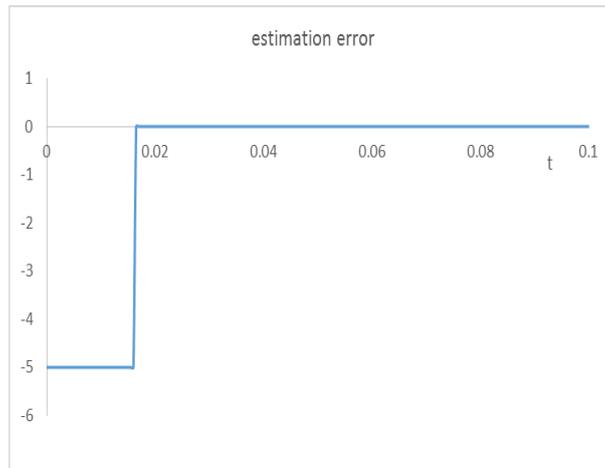

**Fig. 5:** The estimation error $\hat{\theta}(t) - \theta$ for the initial transient period of the closed-loop system (6.1) with (6.5), (6.6), (6.12), (6.13), $A_1 = A_2 = 0$.

The event-triggered adaptive control scheme for $a(x) = |x|^2/20$ and certain $\tilde{N} > 1$, $T > 0$ is given by

$$u(t) = k(\hat{\theta}(\tau_i), x(t)) \quad , \quad t \in [\tau_i, \tau_{i+1}), i \in Z_+$$
$$\hat{\theta}(t) = \hat{\theta}(\tau_i) \quad , \quad t \in [\tau_i, \tau_{i+1}), i \in Z_+ \quad (6.5)$$

where $k(\theta, x)$ is given by (6.2), the time sequence is given by

$$\tau_{i+1} = \tau_i + \min(T, r_i) \quad , \quad i \in Z_+$$
$$\tau_0 = 0 \quad (6.6)$$

and the event trigger is given by:

$$r_i := \inf\left\{ t > \tau_i : V_{\hat{\theta}(\tau_i)}(x(t)) = V_{\hat{\theta}(\tau_i)}(x(\tau_i)) + \frac{1}{20}|x(\tau_i)|^2 \right\}, \text{ for } x(\tau_i) \neq 0 \quad (6.7)$$

$$r_i := T, \text{ for } x(\tau_i) = 0 \quad (6.8)$$

where $V(\theta, x)$ is given by (6.3). The parameter update law for $i \in Z_+$ is given by:

$$\hat{\theta}(\tau_{i+1}) = \hat{\theta}(\tau_i), \text{ if } \Xi(\tau_{i+1}, \mu_{i+1}) = 0 \quad (6.9)$$

$$\hat{\theta}(\tau_{i+1}) = \frac{\int_{\mu_{i+1}}^{\tau_{i+1}} \int_{\mu_{i+1}}^{\tau_{i+1}} \left(\int_\sigma^t x_1^2(s)\,ds\right)\left(x_1(t) - x_1(\sigma) - \int_\sigma^t x_2(s)\,ds\right) d\sigma\, dt}{\Xi(\tau_{i+1}, \mu_{i+1})}, \text{ if } \Xi(\tau_{i+1}, \mu_{i+1}) > 0 \quad (6.10)$$

where

$$\Xi(\tau_{i+1}, \mu_{i+1}) = \int_{\mu_{i+1}}^{\tau_{i+1}} \int_{\mu_{i+1}}^{\tau_{i+1}} \left(\int_\sigma^t x_1^2(s)\,ds\right)^2 d\sigma\, dt, \quad \mu_{i+1} := \min\left\{\tau_j : j \in \{0,\ldots,i\}, \tau_j \geq \tau_{i+1} - \tilde{N}T\right\} \quad (6.11)$$

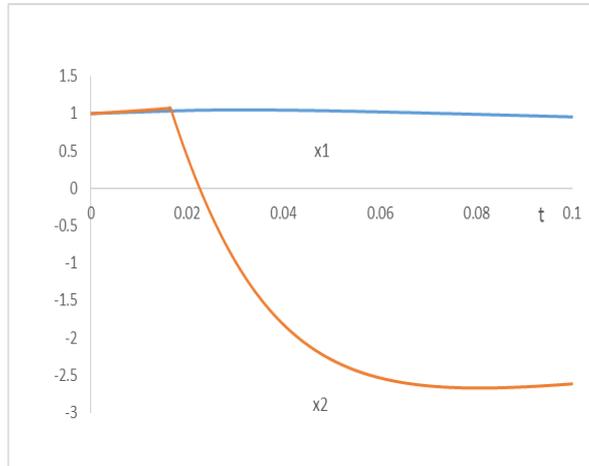

**Fig. 6:** The initial transient period of the closed-loop system (6.1) with (6.5), (6.6), (6.12), (6.13), $A_1 = A_2 = 0$.

Theorem 4.1 and Theorem 4.2 guarantee that:
- there exists a family of $KL$ mappings $\tilde{\sigma}_{\theta,\hat{\theta}} \in KL$ parameterized by $\theta \in \Re$, $\hat{\theta} \in \Re$ such that for every $\theta \in \Re$, $x_0 \in \Re^2$, $\hat{\theta}_0 \in \Re$ the solution of the hybrid closed-loop system (6.1) with $v \equiv 0$, (6.5), (6.6), (6.7), (6.8), (6.9), (6.11) and initial conditions $x(0) = x_0$, $\hat{\theta}(0) = \hat{\theta}_0$ is unique, is defined for all $t \geq 0$ and satisfies $|x(t)| \leq \tilde{\sigma}_{\theta,\hat{\theta}_0}(|x_0|, t)$ for all $t \geq 0$,
- if $x_0 \neq 0$ then $\hat{\theta}(t) = \theta$ for all $t \geq T$.
- there exists a family of constants $\tilde{M}_{\theta,\hat{\theta}}, \omega_\theta, \tilde{R}_{\theta,\hat{\theta}} > 0$ parameterized by $\theta \in \Re$, $\hat{\theta} \in \Re$, such that for every $\theta \in \Re$, $x_0 \in \Re^2$, $\hat{\theta}_0 \in \Re$ with $|x_0| \leq \tilde{R}_{\theta,\hat{\theta}_0}$ the solution of the hybrid closed-loop system (6.1) with $v \equiv 0$, (6.5), (6.6), (6.7), (6.8), (6.9), (6.11) and initial conditions $x(0) = x_0$, $\hat{\theta}(0) = \hat{\theta}_0$ satisfies the estimate $|x(t)| \leq \tilde{M}_{\theta,\hat{\theta}_0} \exp(-\omega_\theta t)|x_0|$ for all $t \geq 0$.

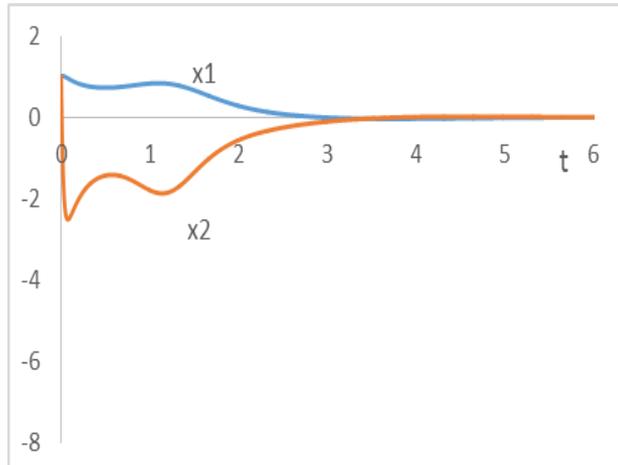

**Fig. 7:** The closed-loop system (6.1) with $u = k(\theta, x)$, $A_1 = 2$, $A_2 = 0$.

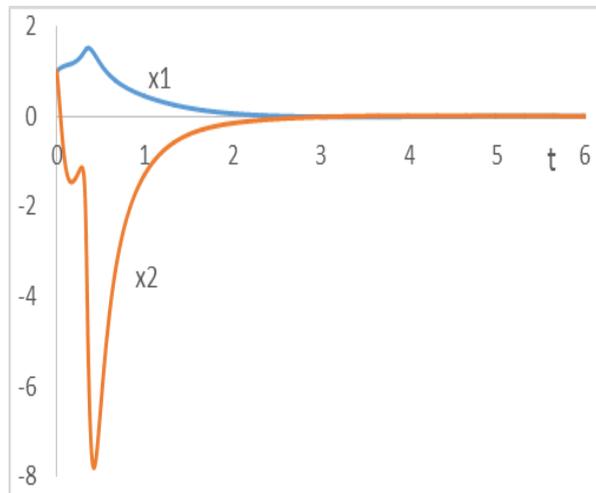

**Fig. 8:** The closed-loop system (6.1) with (6.4), $A_1 = 2$, $A_2 = 0$.

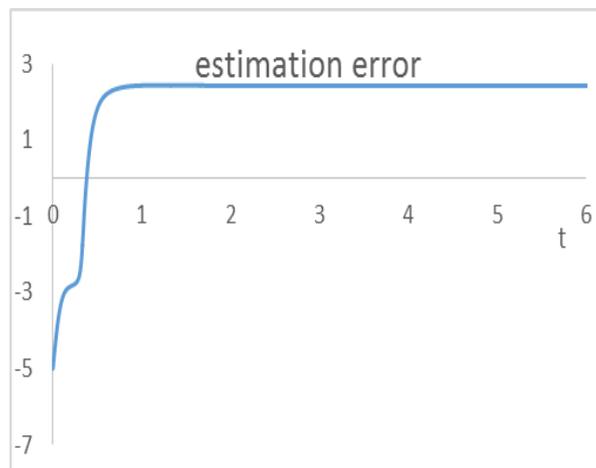

**Fig. 9:** The estimation error $\hat{\theta}(t) - \theta$ for the closed-loop system (6.1) with (6.4), $A_1 = 2$, $A_2 = 0$.

In order to be able to implement numerically the event trigger and the parameter update law, we introduce a small constant $\varepsilon > 0$ and we modify the event trigger and parameter update law (for the initial period $t \leq \tilde{N}T$ and thus $\mu_{i+1} = 0$; see (6.11)) as follows:

$$r_i := \inf\left\{ t > \tau_i : V_{\hat{\theta}(\tau_i)}(x(t)) = V_{\hat{\theta}(\tau_i)}(x(\tau_i)) + \frac{1}{20}|x(\tau_i)|^2 + \varepsilon \right\}, \text{ for all } x(\tau_i) \in \Re^n \quad (6.12)$$

$$\hat{\theta}(\tau_{i+1}) = \begin{cases} \hat{\theta}(\tau_i) & \text{if } \eta(\tau_{i+1}) < \varepsilon \\ \dfrac{\zeta(\tau_{i+1})}{\eta(\tau_{i+1})} & \text{if } \eta(\tau_{i+1}) \geq \varepsilon \end{cases}, \quad i \in Z_+$$

$$\dot{\eta} = 2z_2 + 2z_6 z_1^2 - 4z_1 z_3$$
$$\dot{\zeta} = 2(x_1 - z_7)(z_6 z_1 - z_3) + 2z_5 - 2z_1 z_4 \quad (6.13)$$
$$\dot{z}_1 = x_1^2, \quad \dot{z}_2 = z_1^2, \quad \dot{z}_3 = z_1, \quad \dot{z}_4 = x_1 - z_7$$
$$\dot{z}_5 = (x_1 - z_7)z_1, \quad \dot{z}_6 = 1, \quad \dot{z}_7 = x_2$$

with initial conditions $\eta(0) = \zeta(0) = z_i(0) = 0$ ($i = 1,...,7$). Notice that the differential equations in (6.13) are nothing else but an implementation of the double integrals in (6.10).

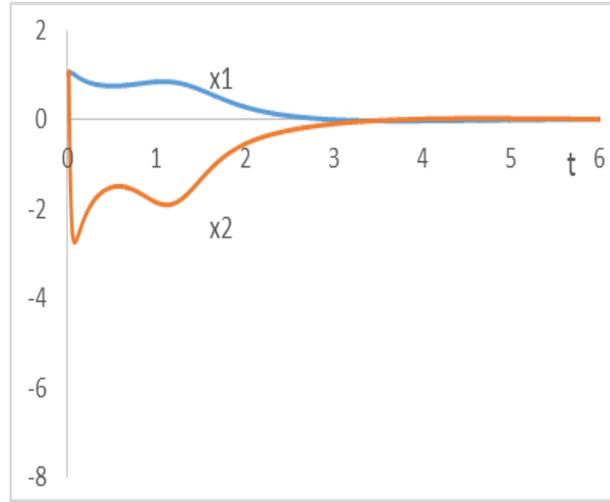

**Fig. 10:** The closed-loop system (6.1) with (6.5), (6.6), (6.12), (6.13), $A_1 = 2$, $A_2 = 0$.

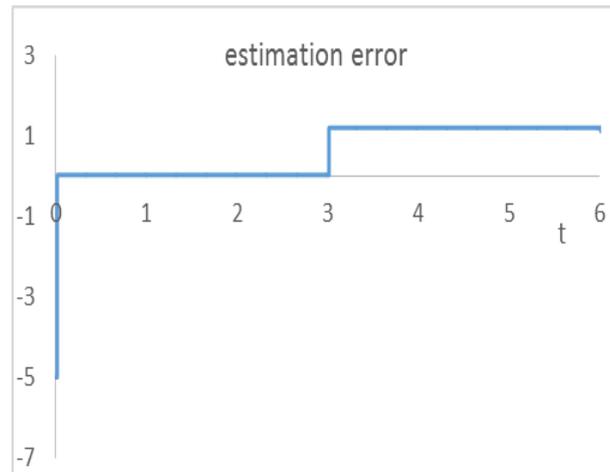

**Fig. 11:** The estimation error $\hat{\theta}(t) - \theta$ for the closed-loop system (6.1) with (6.5), (6.6), (6.12), (6.13), $A_1 = 2$, $A_2 = 0$.

In the following robustness tests we have used $\theta=1$, $\tilde{N}=7$, $\varepsilon=10^{-6}$, $\gamma=5$, $T=3$, $\hat{\theta}(0)=-4$, $x_1(0)=x_2(0)=1$ and

$$v_1(t)=A_1\sin(2t), \quad v_2(t)=A_2\sin(2t). \tag{6.14}$$

First, we consider the disturbance-free case $A_1=A_2=0$. Figure 1 shows the evolution of the state variables for the closed-loop system (6.1) with the nominal feedback law $u=k(\theta,x)$ and known parameter $\theta \in \Re$. On the other hand, Figure 4 shows the evolution of the state variables for the closed-loop system (6.1) with (6.5), (6.6), (6.12), (6.13). The graphs look almost identical. However, they are not identical. In order to see the difference of the responses of the two closed-loop systems we have to check what happens in the initial transient period $t \in [0,0.1]$ for the closed-loop system (6.1) with (6.5), (6.6), (6.12), (6.13): this is shown in Figure 5 and Figure 6. The event trigger is activated at a time close to 0.02 and the parameter $\theta \in \Re$ is estimated exactly. Later than this time, the closed-loop system follows the trajectories of the nominal closed-loop system (6.1) with $u=k(\theta,x)$ and known parameter $\theta \in \Re$. That explains the similarity of Figure 1 and Figure 4. Figure 4 and Figure 5 should also be compared with Figure 2 and Figure 3, where the response of the closed-loop system (6.1) with the adaptive controller (6.4) is shown. It is clear that the estimation of the parameter $\theta \in \Re$ is much slower.

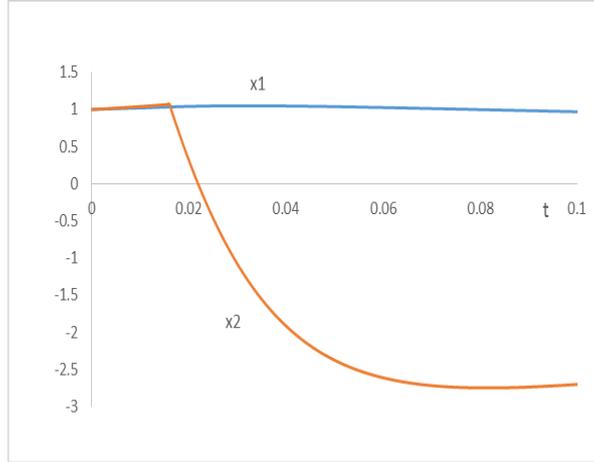

**Fig. 12:** The initial transient period of the closed-loop system (6.1) with (6.5), (6.6), (6.12), (6.13), $A_1=2$, $A_2=0$.

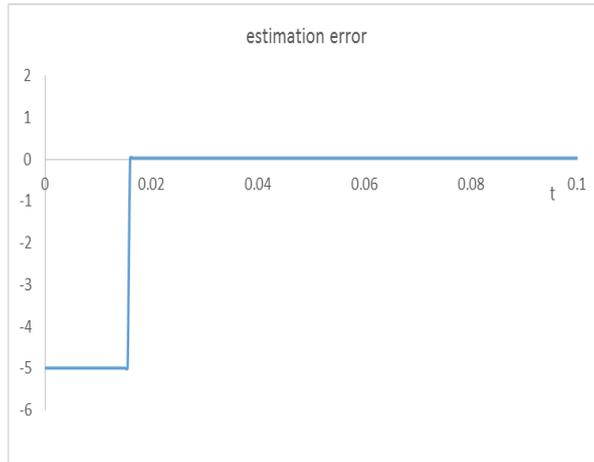

**Fig. 13:** The estimation error $\hat{\theta}(t)-\theta$ for the initial transient period of the closed-loop system (6.1) with (6.5), (6.6), (6.12), (6.13), $A_1=2$, $A_2=0$.

Next, we consider the case of vanishing perturbation $A_1 = 2$, $A_2 = 0$. Figure 7 shows the evolution of the state variables for the closed-loop system (6.1) with the nominal feedback law $u = k(\theta, x)$ and known parameter $\theta \in \Re$. On the other hand, Figure 10 shows the evolution of the state variables for the closed-loop system (6.1) with (6.5), (6.6), (6.12), (6.13). Again, the graphs look almost identical but they are not. In order to see the difference of the responses of the two closed-loop systems we have again to check what happens in the initial transient period $t \in [0, 0.1]$ for the closed-loop system (6.1) with (6.5), (6.6), (6.12), (6.13): this is shown in Figure 12 and Figure 13. The event trigger is activated at a time close to 0.02 and the parameter $\theta \in \Re$ is estimated almost exactly (due to the small effect of the disturbance in this interval). Later than this time and up to time $t = 3$, the closed-loop system follows the trajectories of the nominal closed-loop system (6.1) with $u = k(\theta, x)$ and known parameter $\theta \in \Re$. That explains the similarity of Figure 7 and Figure 10. At a time close to $t = 3$ the event trigger is activated again: see Figure 11. This time the effect of the disturbance is not negligible and leads to an non-negligible estimation error. However, the state has reached a neighborhood around the origin and the effect of the estimation error to the closed-loop system is small. Figure 10 and Figure 11 should also be compared with Figure 8 and Figure 9, where the response of the closed-loop system (6.1) with the adaptive controller (6.4) is shown. It is clear that in this case the adaptive controller (6.4) gives a larger overshoot and a larger estimation error of the parameter.

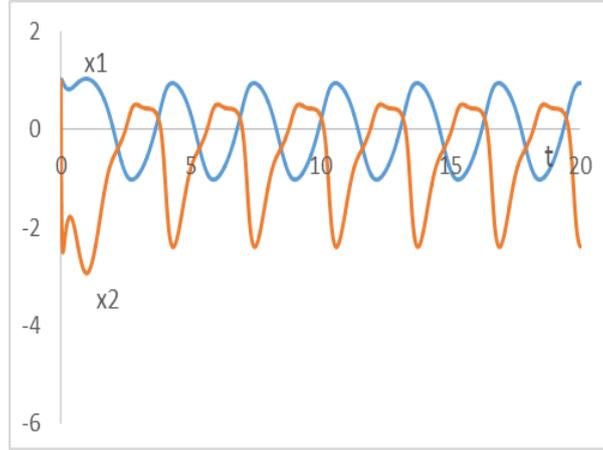

**Fig. 14:** The closed-loop system (6.1) with $u = k(\theta, x)$, $A_2 = 2$, $A_1 = 0$.

Finally, we consider the more demanding case of combined vanishing and non-vanishing perturbation $A_2 = 2$, $A_1 = 0$. Figure 14 shows the evolution of the state variables for the closed-loop system (6.1) with the nominal feedback law $u = k(\theta, x)$ and known parameter $\theta \in \Re$. On the other hand, Figure 17 shows the evolution of the state variables for the closed-loop system (6.1) with (6.5), (6.6), (6.12), (6.13). Again, the graphs look almost identical but they are not. In order to see the difference of the responses of the two closed-loop systems we have again to check what happens in the initial transient period $t \in [0, 0.1]$ for the closed-loop system (6.1) with (6.5), (6.6), (6.12), (6.13): this is shown in Figure 19 and Figure 20. The event trigger is activated at a time close to 0.02 and the parameter is estimated almost exactly (due to the small effect of the disturbances in this interval). Later than this time, the closed-loop system follows the trajectories of the nominal closed-loop system (6.1) with $u = k(\theta, x)$ and known parameter $\theta \in \Re$. That explains the similarity of Figure 14 and Figure 17. The event trigger is activated 35 times in the interval $[0, 20]$: see Figure 18. However, the effect of the disturbances is cancelled, leading to eventually negligible estimation errors. Figure 17 and Figure 18 should also be compared with Figure 15 and Figure 16, where the response of the closed-loop system (6.1) with the adaptive controller (6.4) is shown. It is clear that the adaptive controller (6.4) gives a larger overshoot and eventually leads the system to a limit cycle with a larger oscillation magnitude than the event-triggered adaptive controller (6.5), (6.6), (6.12), (6.13). Moreover, the adaptive controller (6.4) leads to a larger estimation error of the parameter than the event-triggered adaptive controller (6.5), (6.6), (6.12), (6.13).

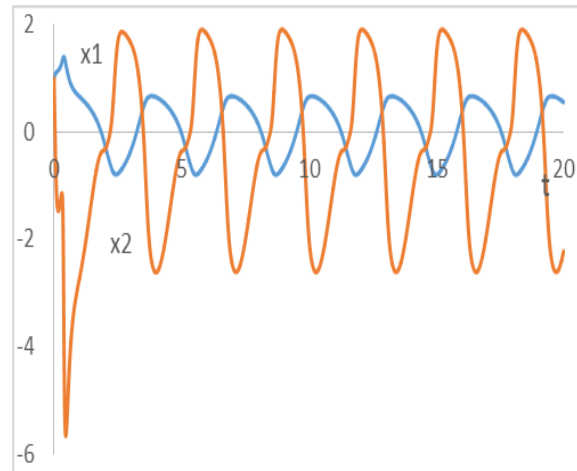

**Fig. 15:** The closed-loop system (6.1) with (6.4), $A_2 = 2$, $A_1 = 0$.

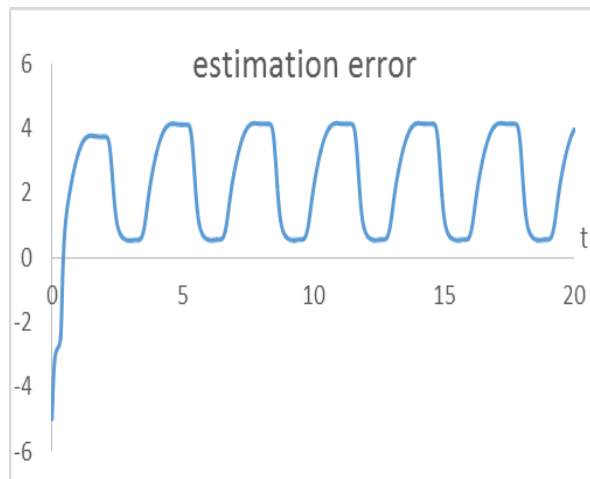

**Fig. 16:** The estimation error $\hat{\theta}(t) - \theta$ for the closed-loop system (6.1) with (6.4), $A_2 = 2$, $A_1 = 0$.

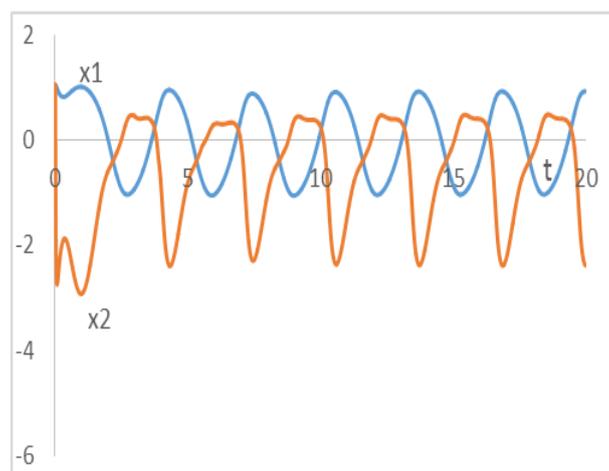

**Fig. 17:** The closed-loop system (6.1) with (6.5), (6.6), (6.12), (6.13), $A_2 = 2$, $A_1 = 0$.

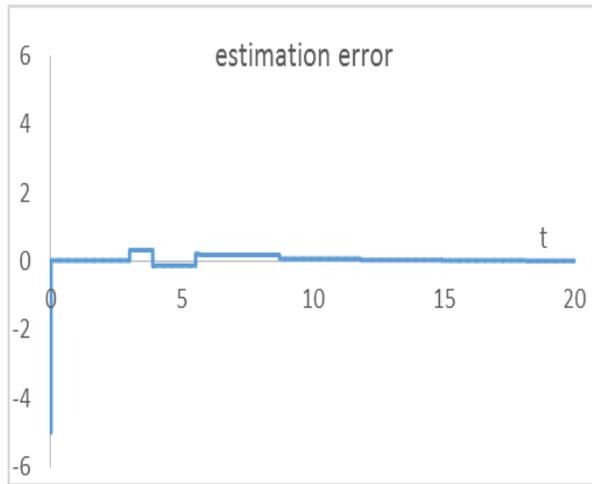

**Fig. 18:** The estimation error $\hat{\theta}(t) - \theta$ for the closed-loop system (6.1) with (6.5), (6.6), (6.12), (6.13), $A_2 = 2$, $A_1 = 0$.

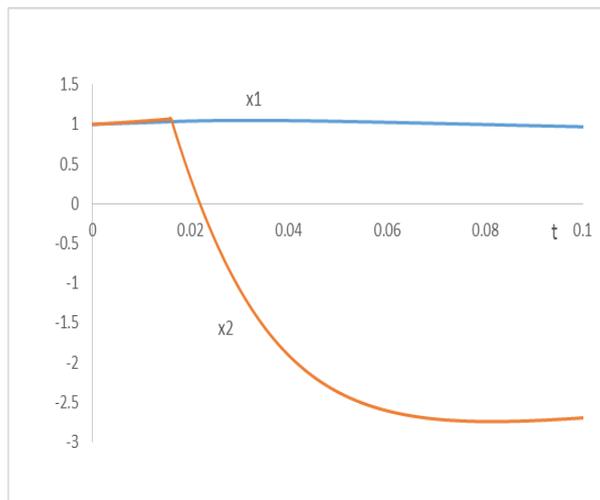

**Fig. 19:** The initial transient period of the closed-loop system (6.1) with (6.5), (6.6), (6.12), (6.13), $A_2 = 2$, $A_1 = 0$.

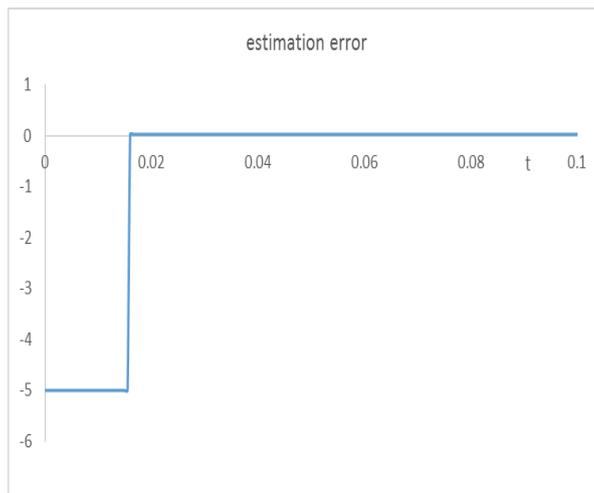

**Fig. 20:** The estimation error $\hat{\theta}(t) - \theta$ for the initial transient period of the closed-loop system (6.1) with (6.5), (6.6), (6.12), (6.13), $A_2 = 2$, $A_1 = 0$.

The conclusions of this simulation study are:
- The proposed event-triggered adaptive controller gives responses which follow closely the trajectories of the closed-loop system with the nominal controller and known parameter values. The responses are better than the corresponding responses of conventional adaptive schemes. The robustness properties of the event-triggered adaptive controller are comparable to the robustness properties of the closed-loop system with the nominal controller and known parameter values.

- The proposed event-triggered adaptive controller gives an almost exact estimation of the parameters even in cases where serious perturbations are active. Moreover, in all cases the estimation of the parameters is performed much faster and much more accurately than in conventional adaptive schemes.

## 7. Concluding Remarks

Our certainty equivalence (indirect) adaptive design with finite-time LS-based regulation-triggered identification is distinct from both the existing tuning-based approaches (Lyapunov-based and modular) and from the supervisory (direct, non-tuning-based) approaches. The major distinctions from the latter are in our indirect approach and the triggering that is based on the control objective (the regulation error). The major distinction from the former is that, for the first time, a regulation error-based parameter tuning is performed outside of a Lyapunov-based context and, more specifically, without requiring a complex redesign of the controller in a non-modular manner where the controller incorporates the update law and compensates for the detrimental effects of rapid tuning on the transients. This is the first use of a simple CE adaptive controller for general nonlinear systems, without restrictions on the growth of nonlinearities in the system's model.

Our dead-beat regulation-triggered LS-based identifier achieves parameter convergence, in finite time, for all initial conditions other than those where the initial plant state is zero, even without PE, and the adaptive controller guarantees global regulation of the plant state. This is likely the first such result in adaptive control literature.

The companion paper [24] provides the proofs of all results as well as a convenient algorithmic way of checking the parameter-observability assumption (H3) for a certain class of nonlinear control systems.

In future work we will address robustness beyond the intuitive dead zone-like modification in (6.12) and the simulation study for vanishing and non-vanishing perturbations in Section 7. Moreover, we will extend the approach to the case of output feedback controllers and we will study adaptive observers.

**Acknowledgments:** The authors would like to thank Professor John Tsinias for the useful exchange of ideas and for suggesting Corollary 4.4.